# Markov Chain Monte Carlo: Can We Trust the Third Significant Figure?

**James M. Flegal, Murali Haran and Galin L. Jones**

*Abstract.* Current reporting of results based on Markov chain Monte Carlo computations could be improved. In particular, a measure of the accuracy of the resulting estimates is rarely reported. Thus we have little ability to objectively assess the quality of the reported estimates. We address this issue in that we discuss why Monte Carlo standard errors are important, how they can be easily calculated in Markov chain Monte Carlo and how they can be used to decide when to stop the simulation. We compare their use to a popular alternative in the context of two examples.

## 1. INTRODUCTION

Hoaglin and Andrews (1975) consider the general problem of what information should be included in publishing computation-based results. The goal of their suggestions was "...to make it easy for the reader to make reasonable assessments of the numerical quality of the results." In particular, Hoaglin and Andrews suggested that it is crucial to report some notion of the accuracy of the results and, for Monte Carlo studies this should include estimated standard errors. However, in settings where Markov chain Monte Carlo (MCMC) is used there is a culture of rarely reporting such information. For example, we looked at the issues published in 2006 of *Journal of the American Statistical Association*, *Biometrika* and *Journal of the Royal Statistical Society, Series B*. In these journals we found 39 papers that used MCMC. Only three of them directly addressed the Monte Carlo error in the reported estimates. Thus it is apparent that the readers of the other papers have little ability to objectively assess the quality of the reported estimates. We attempt to address this issue in that we discuss why Monte Carlo standard errors are important, how they can be easily calculated in MCMC settings and compare their use to a popular alternative.

Simply put, MCMC is a method for using a computer to generate data and subsequently using standard large sample statistical methods to estimate fixed, unknown quantities of a given target distribution. (Thus, we object to calling it "Bayesian Computation.") That is, it is used to produce a point estimate of some characteristic of a target distribution $\pi$ having support $\mathsf{X}$. The most common use of MCMC is to estimate $E_\pi g := \int_\mathsf{X} g(x) \pi(dx)$ where $g$ is a real-valued, $\pi$-integrable function on $\mathsf{X}$.

Suppose that $X = \{X_1, X_2, X_3, \ldots\}$ is an aperiodic, irreducible, positive Harris recurrent Markov chain with state space $\mathsf{X}$ and invariant distribution $\pi$ (for definitions see Meyn and Tweedie (1993)). In this case $X$ is *Harris ergodic*. Typically, estimating $E_\pi g$ is natural since an appeal to the Ergodic Theorem implies that if $E_\pi |g| < \infty$, then, with probabil-ity 1,

$$\bar{g}_n := \frac{1}{n}\sum_{i=1}^{n} g(X_i) \to E_\pi g \quad \text{as } n \to \infty. \quad (1)$$

*James Flegal is Assistant Professor, Department of Statistics, University of California, Riverside, California 92521, USA e-mail: jflegal@ucr.edu. Murali Haran is Assistant Professor, Department of Statistics, The Pennsylvania State University, University Park, Pennsylvania 16802, USA e-mail: mharan@stat.psu.edu. Galin Jones is Associate Professor, School of Statistics, University of Minnesota, Minneapolis, Minnesota 55455, USA e-mail: galin@stat.umn.edu.*

The MCMC method entails constructing a Markov chain $X$ satisfying the regularity conditions described above and then simulating $X$ for a finite number of steps, say $n$, and using $\bar{g}_n$ to estimate $E_\pi g$. The popularity of MCMC largely is due to the ease with which such an $X$ can be simulated (Chen, Shao and Ibrahim (2000), Robert and Casella (1999), Liu (2001)).

An obvious question is when should we stop the simulation? That is, how large should $n$ be? Or, when is $\bar{g}_n$ a good estimate of $E_\pi g$? In a given application we usually have an idea about how many significant figures we want in our estimate, but how should this be assessed? Responsible statisticians and scientists want to do the right thing but output analysis in MCMC has become a muddled area with often conflicting advice and dubious terminology. This leaves many in a position where they feel forced to rely on intuition, folklore and heuristics. We believe this often leads to some poor practices: (A) stopping the simulation too early, (B) wasting potentially useful samples, and, most importantly, (C) providing no notion of the quality of $\bar{g}_n$ as an estimate of $E_\pi g$. In this paper we focus on issue (C) but touch briefly on (A) and (B).

The rest of this paper is organized as follows. In Section 2 we briefly introduce some basic concepts from the theory of Markov chains. In Section 3 we consider estimating the Monte Carlo error of $\bar{g}_n$. Then Section 4 covers two methods for stopping the simulation and compares them in a toy example. In Section 5 the two methods are compared again in a realistic spatial model for a data set on wheat crop flowering dates in North Dakota. We close with some final remarks in Section 6.

## 2. MARKOV CHAIN BASICS

Suppose that $X = \{X_1, X_2, \ldots\}$ is a Harris ergodic Markov chain with state space $\mathsf{X}$ and invariant distribution $\pi$. For $n \in \mathbb{N} := \{1, 2, 3, \ldots\}$ let $P^n(x, \cdot)$ be the $n$-step Markov transition kernel; that is, for $x \in \mathsf{X}$ and a measurable set $A$, $P^n(x, A) = \Pr(X_{n+i} \in A \mid X_i = x)$. An extremely useful property of $X$ is that the chain will converge to the invariant distribution. Specifically,

$$\|P^n(x, \cdot) - \pi(\cdot)\| \downarrow 0 \quad \text{as } n \to \infty,$$

where the left-hand side is the *total variation* distance between $P^n(x, \cdot)$ and $\pi(\cdot)$. (This is stronger than convergence in distribution.) The Markov chain $X$ is *geometrically ergodic* if there exists a constant $0 < t < 1$ and a function $M : \mathsf{X} \to \mathbb{R}^+$ such that for any $x \in \mathsf{X}$,

$$\|P^n(x, \cdot) - \pi(\cdot)\| \le M(x) t^n \quad (2)$$

for $n \in \mathbb{N}$. If $M(x)$ is bounded, then $X$ is *uniformly ergodic*. Thus uniform ergodicity implies geometric ergodicity. However, as one might imagine, finding $M$ and $t$ directly is often quite difficult in realistic settings.

There has been a substantial amount of effort devoted to establishing (2) in MCMC settings. For example, Hobert and Geyer (1998), Johnson and Jones (2008), Jones and Hobert (2004), Marchev and Hobert (2004), Mira and Tierney (2002), Robert (1995), Roberts and Polson (1994), Roberts and Rosenthal (1999), Rosenthal (1995, 1996), Roy and Hobert (2007) and Tierney (1994) examined Gibbs samplers while Christensen, Moller and Waagepetersen (2001), Douc et al. (2004), Fort and Moulines (2000, 2003), Geyer (1999), Jarner and Hansen (2000), Jarner and Roberts (2002), Meyn and Tweedie (1994) and Mengersen and Tweedie (1996) considered Metropolis–Hastings algorithms.

## 3. MONTE CARLO ERROR

A Monte Carlo approximation is not exact. The number $\bar{g}_n$ is not the exact value of the integral we are trying to approximate. It is off by some amount, the *Monte Carlo error*, $\bar{g}_n - E_\pi g$. How large is the Monte Carlo error? Unfortunately, we can never know unless we know $E_\pi g$.

We do not know the Monte Carlo error, but we can get a handle on its sampling distribution. That is, assessing the Monte Carlo error can be accomplished by estimating the variance of the asymptotic distribution of $\bar{g}_n$. Under regularity conditions, the Markov chain $X$ and function $g$ will admit a CLT. That is,

$$\sqrt{n}(\bar{g}_n - E_\pi g) \xrightarrow{d} \mathrm{N}(0, \sigma_g^2) \quad (3)$$

as $n \to \infty$ where $\sigma_g^2 := \mathrm{var}_\pi\{g(X_1)\} + 2\sum_{i=2}^{\infty} \mathrm{cov}_\pi\{g(X_1), g(X_i)\}$; the subscript $\pi$ means that the expectations are calculated assuming $X_1 \sim \pi$. The CLT holds for *any* initial distribution when either (i) $X$ is geometrically ergodic and $E_\pi|g|^{2+\delta} < \infty$ for some $\delta > 0$ or (ii) $X$ is uniformly ergodic and $E_\pi g^2 < \infty$. These are not the only sufficient conditions for a CLT but are among the most straightforward to state; the interested reader is pointed to the summaries provided by Jones (2004) and Roberts and Rosenthal (2004).

Given a CLT we can assess the Monte Carlo error in $\bar{g}_n$ by estimating the variance, $\sigma_g^2$. That is, we can calculate and report an estimate of $\sigma_g^2$, say $\hat{\sigma}_g^2$, that will allow us to assess the accuracy of the point estimate. There have been many techniques introduced for estimating $\sigma_g^2$; see, among others, Bratley, Fox and Schrage (1987), Fishman (1996), Geyer (1992), Glynn and Iglehart (1990), Glynn and Whitt (1991), Mykland, Tierney and Yu (1995) and Roberts (1996). For example, regenerative simulation, batch means and spectral variance estimators all can be appropriate in MCMC settings. We will consider only one of the available methods, namely non-overlapping batch means. We chose this method because it is easy to implement and can enjoy desirable theoretical properties. However, overlapping batch means has a reputation of sometimes being more efficient than nonoverlapping batch means. On the other hand, currently the spectral variance and overlapping batch means estimators require stronger regularity conditions than nonoverlapping batch means.

### 3.1 Batch Means

In nonoverlapping batch means the output is broken into blocks of equal size. Suppose the algorithm is run for a total of $n = ab$ iterations (hence $a = a_n$ and $b = b_n$ are implicit functions of $n$) and define

$$\bar{Y}_j := \frac{1}{b} \sum_{i=(j-1)b+1}^{jb} g(X_i) \quad \text{for } j = 1, \ldots, a.$$

The batch means estimate of $\sigma_g^2$ is

$$(4) \quad \hat{\sigma}_g^2 = \frac{b}{a-1} \sum_{j=1}^{a} (\bar{Y}_j - \bar{g}_n)^2.$$

Batch means is attractive because it is easy to implement (and it is available in some software, e.g., WinBUGS) but some authors encourage caution in its use (Roberts (1996)). In particular, we believe careful use is warranted since (4), in general, is not a consistent estimator of $\sigma_g^2$. On the other hand, Jones et al. (2006) showed that if the batch size and the number of batches are allowed to increase as the overall length of the simulation increases by setting $b_n = \lfloor n^\theta \rfloor$ and $a_n = \lfloor n/b_n \rfloor$, then $\hat{\sigma}_g^2 \to \sigma_g^2$ with probability 1 as $n \to \infty$. In this case we call it consistent batch means (CBM) to distinguish it from the standard (fixed number of batches) version. The regularity conditions require that $X$ be geometrically ergodic, $E_\pi |g|^{2+\varepsilon_1+\varepsilon_2} < \infty$ for some $\varepsilon_1 > 0$, $\varepsilon_2 > 0$ and $(1 + \varepsilon_1/2)^{-1} < \theta < 1$; often $\theta = 1/2$ (i.e., $b_n = \lfloor \sqrt{n} \rfloor$ and $a_n = \lfloor n/b_n \rfloor$) is a convenient choice that works well in applications. Note that the only practical difference between CBM and standard batch means is that the batch number and size are chosen as functions of the overall run length, $n$. A simple R function for implementing CBM or a faster command line C version of this function is available from the authors upon request.

Using CBM to get an estimate of the Monte Carlo standard error (MCSE) of $\bar{g}_n$, say $\hat{\sigma}_g/\sqrt{n}$, we can form an asymptotically valid confidence interval for $E_\pi g$. The half-width of the interval is given by

$$(5) \quad t_{a_n-1} \frac{\hat{\sigma}_g}{\sqrt{n}}$$

where $t_{a_n-1}$ is an appropriate quantile from Student's $t$ distribution with $a_n - 1$ degrees of freedom.

### 3.2 How Many Significant Figures?

The title of the paper contains a rhetorical question; we do not always care about the *third* significant figure. But we should care about how many significant figures there are in our estimates. Assessing the Monte Carlo error through (5) gives us a tool to do this. For example, suppose $\bar{g}_n = 0.02$; then there is exactly one significant figure in the estimate, namely the "2," but how confident are we about it? Letting $h_\alpha$ denote the half-width given in (5) of a $(1-\alpha)100\%$ interval, we would trust the one significant figure in our estimate if $0.02 \pm h_\alpha \subseteq [0.015, 0.025)$ since otherwise values such as $E_\pi g = 0.01$ or $E_\pi g = 0.03$ are plausible through rounding.

More generally, we can use (5) to assess how many significant figures we have in our estimates. This is illustrated in the following toy example that will be used several times throughout the rest of this paper.

3.2.1 *Toy example.* Let $Y_1, \ldots, Y_K$ be i.i.d. $N(\mu, \lambda)$ and let the prior for $(\mu, \lambda)$ be proportional to $1/\sqrt{\lambda}$. The posterior density is characterized by

$$(6) \quad \pi(\mu, \lambda | y) \propto \lambda^{-(K+1)/2} \exp\left\{-\frac{1}{2\lambda} \sum_{j=1}^{K} (y_j - \mu)^2\right\}$$

where $y = (y_1, \ldots, y_K)^T$. It is easy to check that this posterior is proper as long as $K \geq 3$ and we assume this throughout. Using the Gibbs sampler to make draws from (6) requires the full conditional densities, $f(\mu|\lambda, y)$ and $f(\lambda|\mu, y)$, which are as follows:

$$\mu|\lambda, y \sim \text{N}(\bar{y}, \lambda/K),$$
$$\lambda|\mu, y \sim \text{IG}\left(\frac{K-1}{2}, \frac{(K-1)s^2 + K(\bar{y}-\mu)^2}{2}\right),$$

where $\bar{y}$ is the sample mean and $(K-1)s^2 = \sum(y_i - \bar{y})^2$. [We say $W \sim \text{IG}(\alpha, \beta)$ if its density is proportional to $w^{-(\alpha+1)}e^{-\beta/w}I(w > 0)$.] Consider estimating the posterior means of $\mu$ and $\lambda$. It is easy to prove that $E(\mu|y) = \bar{y}$ and $E(\lambda|y) = (K-1)s^2/(K-4)$ for $K > 4$. Thus we do not need MCMC to estimate these quantities but we will ignore this and use the output of a Gibbs sampler to estimate $E(\mu|y)$ and $E(\lambda|y)$.

Consider the Gibbs sampler that updates $\lambda$ then $\mu$; that is, letting $(\lambda', \mu')$ denote the current state and $(\lambda, \mu)$ denote the future state, the transition looks like $(\lambda', \mu') \to (\lambda, \mu') \to (\lambda, \mu)$. Jones and Hobert (2001) established that the associated Markov chain is geometrically ergodic as long as $K \geq 5$. If $K > 10$, then the moment conditions ensuring the CLT and the regularity conditions for CBM (with $\theta = 1/2$) hold.

Let $K = 11$, $\bar{y} = 1$, and $(K-1)s^2 = 14$ so that $E(\mu|y) = 1$ and $E(\lambda|y) = 2$; for the remainder of this paper these settings will be used every time we consider this example. Consider estimating $E(\mu|y)$ and $E(\lambda|y)$ with $\bar{\mu}_n$ and $\bar{\lambda}_n$, respectively, and using CBM to calculate the MCSEs for each estimate. Specifically, we will use a 95% confidence level in (5) to construct an interval estimate. Let the initial value for the simulation be $(\lambda_1, \mu_1) = (1, 1)$. When we ran the Gibbs sampler for 1000 iterations we obtained $\bar{\lambda}_{1000} = 2.003$ with an MCSE of 0.055 and $\bar{\mu}_{1000} = 0.99$ with an MCSE of 0.016. Thus we would be comfortable reporting two significant figures for the estimates of $E(\lambda|y)$ and $E(\mu|y)$, specifically 2.0 and 1.0, respectively. But when we started from $(\lambda_1, \mu_1) = (100, 100)$ and ran Gibbs for 1000 iterations we obtained $\bar{\lambda}_{1000} = 13.06$ with an MCSE of 11.01 and $\bar{\mu}_{1000} = 1.06$ with an MCSE of 0.071. Thus we would not be comfortable with *any* significant figures for the estimate of $E(\lambda|y)$ but we would trust one significant figure (i.e., 1) for $E(\mu|y)$. Unless the MCSE is calculated and reported a hypothetical reader would have no way to judge this independently.

3.2.2 *Remarks.*

1. A common concern about MCSEs is that their use may require estimating $E_\pi g$ much too precisely relative to $\sqrt{\text{var}_\pi g}$. Of course, it would be a rare problem indeed where we would know $\sqrt{\text{var}_\pi g}$ and not $E_\pi g$. Thus we would need to estimate $\sqrt{\text{var}_\pi g}$ and calculate an MCSE (via the delta method) before we could trust the estimate of $\sqrt{\text{var}_\pi g}$ to inform us about the MCSE for $E_\pi g$.
2. We are not suggesting that all MCMC-based estimates should be reported in terms of significant figures; in fact we do not do this in the simulations that occur later. Instead, we are strongly suggesting that an estimate of the Monte Carlo standard error should be used to assess simulation error and reported. Without an attached MCSE a point estimate should not be trusted.

## 4. STOPPING THE SIMULATION

In this section we consider two formal approaches to terminating the simulation. The first is based on calculating an MCSE and is discussed in Section 4.1. The second is based on the method introduced in Gelman and Rubin (1992) and is one of many so-called convergence diagnostics (Cowles and Carlin (1996)). Our reason for choosing the Gelman–Rubin diagnostic (GRD) is that it appears to be far and away the most popular method for stopping the simulation. GRD and MCSE are used to stop the simulation in a similar manner. After $n$ iterations either the value of the GRD or MCSE is calculated and if it is not sufficiently small then we continue the simulation until it is.

### 4.1 Fixed-Width Methodology

Suppose we have an idea of how many significant figures we want in our estimate. Another way of saying this is that we want the half-width of the interval (5) to be less than some user-specified value, $\varepsilon$. Thus we might consider stopping the simulation when the MCSE of $\bar{g}_n$ is sufficiently small. This, of course, means that we may have to check whether this criterion is met many times. It is not obvious that such a procedure will be guaranteed to terminate the computation in a finite amount of time or whether the resulting intervals will enjoy the desired coverage probability and half-width. Also, we do not want to check too early in the simulation since we will run the risk of premature termination due to a poor estimate of the standard error.

Suppose we use CBM to estimate the Monte Carlo standard error of $\bar{g}_n$, say $\hat{\sigma}_g/\sqrt{n}$, and use it to form a confidence interval for $E_\pi g$. If this interval is too large, then the value of $n$ is increased and simulation continues until the interval is sufficiently small. Formally, the criterion is given by

$$(7) \qquad t_{a_n-1} \frac{\hat{\sigma}_g}{\sqrt{n}} + p(n) \leq \varepsilon$$

where $t_{a_n-1}$ is an appropriate quantile, $p(n) = \varepsilon I(n < n^*)$ where $n^* > 0$ is fixed, $I$ is the usual indicator function on $\mathbb{Z}_+$ and $\varepsilon > 0$ is the user-specified half-width. The role of $p$ is to ensure that the simulation is not terminated prematurely due to a poor estimate of $\hat{\sigma}_g$. The conditions which guarantee $\hat{\sigma}_g^2$ is consistent also imply that this procedure will terminate in a finite amount of time with probability 1 and that the resulting intervals asymptotically have the desired coverage (see Glynn and Whitt (1992)). However, the finite-sample properties of (5) have received less formal investigation but simulation results suggest that the resulting intervals have approximately the desired coverage and half-width in practice (see Jones et al., 2006).

4.1.1 *Remarks.*

1. The CLT and CBM require a geometrically ergodic Markov chain. This can be difficult to check directly in any given application. On the other hand, considerable effort has been spent establishing (2) for a number of Markov chains; see the references given at the end of Section 2. In our view, this is not the obstacle that it was in the past.
2. The frequency with which (7) should be evaluated is an open question. Checking often, say every few iterations, may substantially increase the overall computational effort.
3. Consider $p(n) = \varepsilon I(n < n^*)$. The choice of $n^*$ is often made based on the user's experience with the problem at hand. However, for geometrically ergodic Markov chains there is some theory that can give guidance on this issue (see Jones and Hobert (2001), Rosenthal (1995)).
4. Stationarity of the Markov chain is not required for the CLT or the strong consistency of CBM. One consequence is that burn-in is not required if we can find a reasonable starting value.

4.1.2 *Toy example.* We consider implementation of fixed-width methods in the toy example introduced in Section 3.2.1. We performed 1000 independent replications of the following procedure. Each replication of the Gibbs sampler was started from $\bar{y}$. Using (7), a replication was terminated when the half-width of a 95% interval with $p(n) = \varepsilon I(n < 400)$ was smaller than a prespecified cutoff, $\varepsilon$, for *both* parameters. If both standard errors were not less than the cutoff, then the current chain length was increased by 10% before checking again. We used two settings for the cutoff, $\varepsilon = 0.06$ and $\varepsilon = 0.04$. These settings will be denoted CBM1 and CBM2, respectively.

First, consider the estimates of $E(\mu|y)$. We can see from Figure 1(a) and (b) that the estimates of $E(\mu|y)$ are centered around the truth with both settings. Clearly, the cutoff of $\varepsilon = 0.04$ is more stringent and yields estimates that are closer to the true value. It should come as no surprise that the cost of this added precision is increased computational effort; see Table 2. The corresponding plots for $\bar{\lambda}_n$ yield the same results and are therefore excluded.

Consider CBM2. In this case, 100% of the estimates, $\bar{\mu}_n$, of $E(\mu|y)$ and 96% of the estimates, $\bar{\lambda}_n$, of $E(\lambda|y)$ are within the specified $\varepsilon = 0.04$ of the truth. In every replication the simulation was stopped when the criterion (7) for $E(\lambda|y)$ dropped below the cutoff. Similar results hold for the CBM1 ($\varepsilon = 0.06$) setting.

## 4.2 The Gelman–Rubin Diagnostic

The Gelman–Rubin diagnostic (GRD) introduced in Gelman and Rubin (1992) and refined by Brooks and Gelman (1998) is a popular method for assessing the output of MCMC algorithms. It is important to note that this method is also based on a Markov chain CLT (Gelman and Rubin (1992), page 463) and hence does not apply more generally than approaches based on calculating an MCSE.

GRD is based on the simulation of $m$ independent parallel Markov chains having invariant distribution $\pi$, each of length $2l$. Thus the total simulation effort is $2lm$. Gelman and Rubin (1992) suggest that the first $l$ simulations should be discarded and inference based on the last $l$ simulations; for the $j$th chain these are denoted $\{X_{1j}, X_{2j}, X_{3j}, \ldots, X_{lj}\}$ with $j = 1, 2, \ldots, m$. Recall that we are interested in estimating $E_\pi g$ and define $Y_{ij} = g(X_{ij})$,

$$B = \frac{l}{m-1} \sum_{j=1}^{m} (\bar{Y}_{\cdot j} - \bar{Y}_{\cdot \cdot})^2 \quad \text{and} \quad W = \frac{1}{m} \sum_{j=1}^{m} s_j^2$$

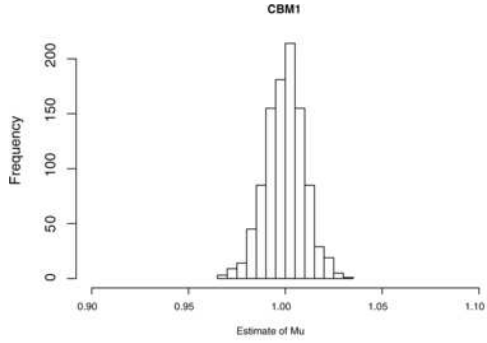

(a) CBM1, with a cutoff of $\epsilon = 0.06$.

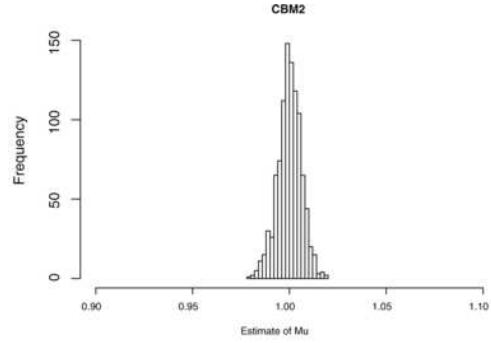

(b) CBM2, with a cutoff of $\epsilon = 0.04$.

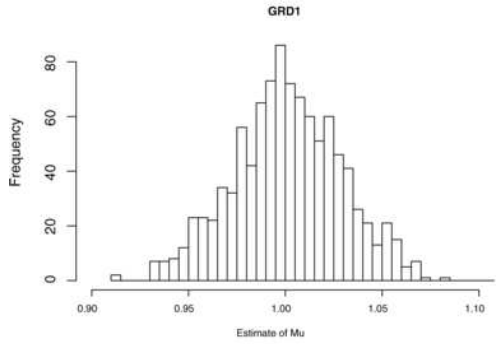

(c) GRD1, 2 chains and $\delta = 1.1$.

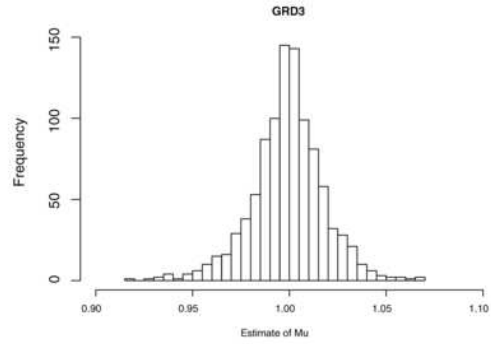

(d) GRD3, 2 chains and $\delta = 1.005$.

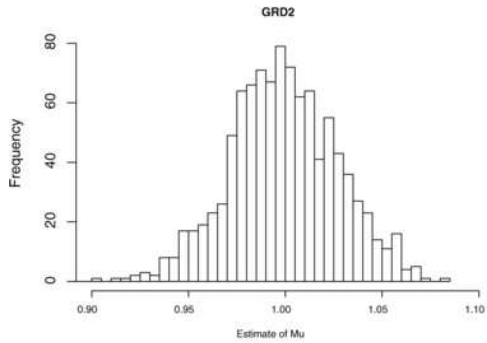

(e) GRD2, 4 chains and $\delta = 1.1$.

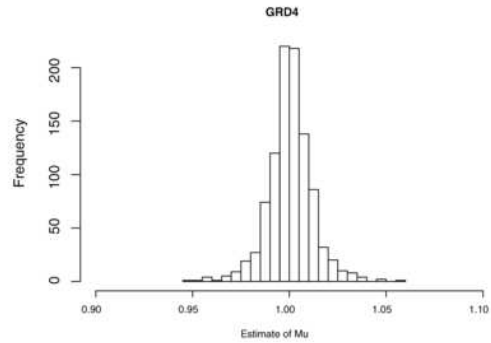

(f) GRD4, 4 chains and $\delta = 1.005$.

FIG. 1. *Histograms from 1000 replications estimating $E(\mu|y)$ for the toy example of Section 3.2.1 with CBM and GRD. Simulation sample sizes are given in Table 2.*

where $\bar{Y}_{.j} = l^{-1}\sum_{i=1}^{l} Y_{ij}$, $\bar{Y}_{..} = m^{-1}\sum_{j=1}^{m} \bar{Y}_{.j}$ and $s_j^2 = (l-1)^{-1}\sum_{i=1}^{l}(Y_{ij} - \bar{Y}_{.j})^2$. Note that $\bar{Y}_{..}$ is the resulting point estimate of $E_\pi g$. Let

$$\hat{V} = \frac{l-1}{l}W + \frac{(m+1)B}{ml}, \quad d \approx \frac{2\hat{V}}{\hat{\text{var}}(\hat{V})},$$

and define the corrected *potential scale reduction factor*

$$\hat{R} = \sqrt{\frac{d+3}{d+1}\frac{\hat{V}}{W}}.$$

As noted by Gelman et al. (2004), $\hat{V}$ and $W$ are essentially two different estimators of $\text{var}_\pi g$; *not* $\sigma_g^2$ from the Markov chain CLT. That is, neither $\hat{V}$ nor $W$ address the sampling variability of $\bar{g}_n$ and hence neither does $\hat{R}$.

In our examples we used the R package coda which reports an upper bound on $\hat{R}$. Specifically, a 97.5% upper bound for $\hat{R}$ is given by

$$\hat{R}_{0.975} = \sqrt{\frac{d+3}{d+1}\left[\frac{l-1}{l} + F_{0.975,m-1,w}\left(\frac{m+1}{ml}\frac{B}{W}\right)\right]},$$

where $F_{0.975,m-1,w}$ is the 97.5% percentile of an $F_w^{m-1}$ distribution, $w = 2W^2/\hat{\sigma}_W^2$ and

$$\hat{\sigma}_W^2 = \frac{1}{m-1}\sum_{j=1}^m (s_j^2 - W)^2.$$

In order to stop the simulation the user provides a cutoff, $\delta > 0$, and simulation continues until

(8) $\qquad \hat{R}_{0.975} + p(n) \leq \delta.$

As with fixed-width methods, the role of $p(n)$ is to ensure that we do not stop the simulation prematurely due to a poor estimate, $\hat{R}_{0.975}$. By requiring a minimum total simulation effort of $n^* = 2lm$ we are effectively setting $p(n) = \delta I(n < n^*)$ where $n$ indexes the total simulation effort.

4.2.1 *Remarks.*

1. A rule of thumb suggested by Gelman et al. (2004) is to set $\delta = 1.1$. These authors also suggest that a value of $\delta$ closer to 1 will be desirable in a "final analysis in a critical problem" but give no further guidance. Since neither $\hat{R}$ nor $\hat{R}_{0.975}$ directly estimates the Monte Carlo error in $\bar{g}_n$ it is unclear to us that $\hat{R} \approx 1$ implies $\bar{g}_n$ is a good estimate of $E_\pi g$.
2. How large should $m$ be? There seem to be few guidelines in the literature except that $m \geq 2$ since otherwise we cannot calculate $B$. Clearly, if $m$ is too large then each chain will be too short to achieve any reasonable expectation of convergence within a given computational effort.
3. The initial values, $X_{j1}$, of the $m$ parallel chains should be drawn from an "overdispersed" distribution. Gelman and Rubin (1992) suggest estimating the modes of $\pi$ and then using a mixture distribution whose components are centered at these modes. Constructing this distribution could be difficult and is often not done in practice (Gelman et al., 2004, page 593).
4. To our knowledge there has been no discussion in the literature about optimal choices of $p(n)$ or $n^*$. In particular, we know of no guidance about how long each of the parallel chains should be simulated before the first time we check that $\hat{R}_{0.975} < \delta$ or how often one should check after that. However, the same theoretical results that could give guidance in item 3 of Section 4.1.1 would apply here as well.
5. GRD was originally introduced simply as a method for determining an appropriate amount of burn-in. However, using diagnostics in this manner may introduce additional bias into the results; see Cowles, Roberts and Rosenthal (1999).

4.2.2 *Toy example.* We consider implementation of GRD in the toy example introduced in Section 3.2.1. The first issue we face is choosing the starting values for each of the $m$ parallel chains. Notice that

$$\pi(\mu,\lambda|y) = g_1(\mu|\lambda)g_2(\lambda)$$

where $g_1(\mu|\lambda)$ is a $N(\bar{y}, \lambda/K)$ density and $g_2(\lambda)$ is an $IG((K-2)/2, (K-1)s^2/2)$ density. Thus we can sequentially sample the exact distribution by first drawing from $g_2(\lambda)$, and then conditionally, draw from $g_1(\mu|\lambda)$. We will use this to obtain starting values for each of the $m$ parallel chains. Thus each of the $m$ parallel Markov chains will be stationary and hence GRD should be at a slight advantage compared to CBM started from $\bar{y}$.

Our goal is to investigate the finite-sample properties of the GRD by considering the estimates of $E(\mu|y)$ and $E(\lambda|y)$ as in Section 4.1.2. To this end, we took multiple chains starting from different draws from the sequential sampler. The multiple chains were run until the total simulation effort was $n^* = 400$ draws; this is the same minimum simulation effort we required of CBM in the previous section. If $\hat{R}_{0.975} < \delta$ for both, the simulation was stopped. Otherwise, 10% of the current chain length was added to each chain before $\hat{R}_{0.975}$ was recalculated. This continued until $\hat{R}_{0.975}$ was below $\delta$ for both. This simulation procedure was repeated independently 1000 times with each replication using the same initial values. We considered four settings using the combinations of $m \in \{2,4\}$ and $\delta \in \{1.005, 1.1\}$. These settings will be denoted GRD1, GRD2, GRD3 and GRD4; see Table 1 for the different settings along with summary statistics that will be considered later.

Upon completion of each replication, the values of $\bar{\mu}_n$ and $\bar{\lambda}_n$ were recorded. Figure 1(c)–(f) show histograms of $\bar{\mu}_n$ for each setting. We can see that all the settings center around the true value of 1, and setting $\delta = 1.005$ provides better estimates. Increasing the number of chains seems to have little impact on the quality of estimation, particularly when $\delta = 1.1$. Histograms of $\bar{\lambda}_n$ for each setting show similar trends.

In the settings we investigated, GRD often terminated the simulations much sooner than CBM. Table 2 shows the percentage of the 1000 replications which were stopped at their minimum ($n^* = 400$) and the percentage with less than 1000 total draws. The data clearly show that premature stopping was common with GRD but uncommon with CBM. This is especially the case for GRD1 and GRD2 where over half the replications used the minimum simulation effort.

Also, the simulation effort for GRD was more variable than that of CBM. In particular, the average simulation effort was comparable for CBM1 and GRD3 and also CBM2 and GRD4 but the standard errors are larger for GRD. Next, Figure 2(a) and (b) show a plot of the estimates, $\bar{\mu}_n$, versus the total number of draws in the chains for CBM2 and GRD4. The graphs clearly show that the total number of draws was more variable using GRD4 than using CBM2. From a practical standpoint, this implies that when using GRD we are likely to run a simulation either too long or too short. Of course, if we run the simulation too long, we will be likely to get a better estimate. But if the simulation is too short, the estimate can be poor.

TABLE 1
*Summary table for all settings and estimated mean-squared error for estimating $E(\mu|y)$ and $E(\lambda|y)$ for the toy example of Section 3.2.1*

| Method | Chains | Stopping rule | MSE for $E(\mu|y)$ | S.E. | MSE for $E(\lambda|y)$ | S.E. |
|---|---|---|---|---|---|---|
| CBM1 | 1 | 0.06 | 9.82e-05 | 4.7e-06 | 1.03e-03 | 4.5e-05 |
| CBM2 | 1 | 0.04 | 3.73e-05 | 1.8e-06 | 3.93e-04 | 1.8e-05 |
| GRD1 | 2 | 1.1 | 7.99e-04 | 3.6e-05 | 8.7e-03 | 4e-04 |
| GRD2 | 4 | 1.1 | 7.79e-04 | 3.7e-05 | 8.21e-03 | 3.6e-04 |
| GRD3 | 2 | 1.005 | 3.49e-04 | 2.1e-05 | 3.68e-03 | 2e-04 |
| GRD4 | 4 | 1.005 | 1.34e-04 | 9.2e-06 | 1.65e-03 | 1.2e-04 |

Standard errors (S.E.) shown for each estimate.

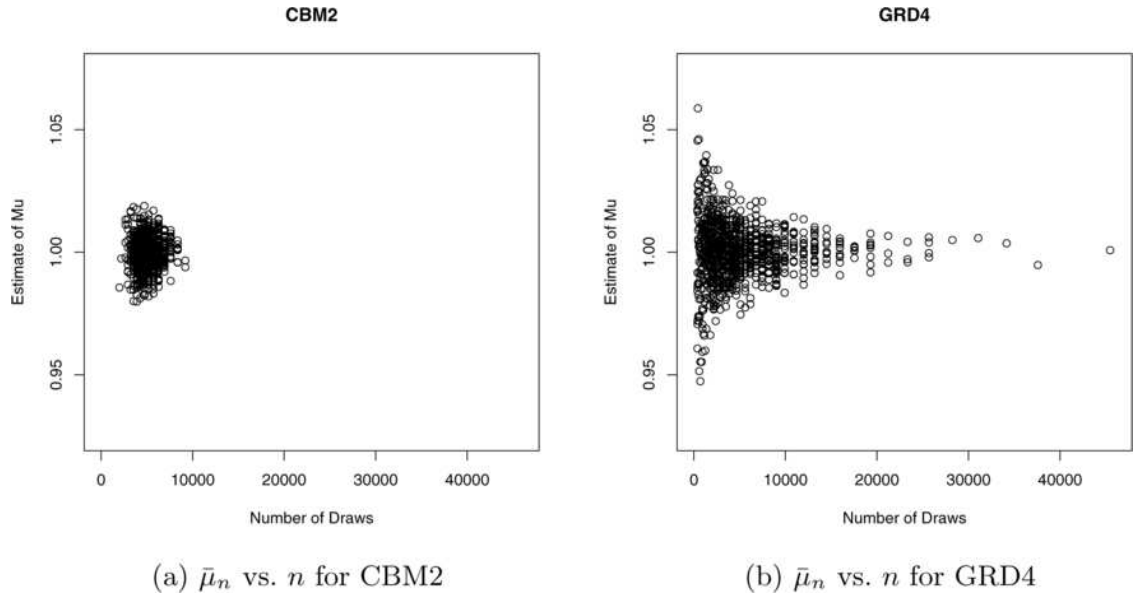

(a) $\bar{\mu}_n$ vs. $n$ for CBM2      (b) $\bar{\mu}_n$ vs. $n$ for GRD4

FIG. 2. *Estimating $E(\mu|y)$ for the toy example of Section 3.2.1. Estimates of $E(\mu|y)$ versus number of draws for the CBM2 and GRD4 settings.*

TABLE 2
*Summary of the proportion (and standard error) of the observed estimates which were based on the minimum number (400) of draws, less than or equal to 1000 draws, and the average total simulation effort for the toy example in Section 3.2.1.*

| Method | Prop. at min. | S.E. | Prop. ≤1000 | S.E. | $N$ | S.E. |
|---|---|---|---|---|---|---|
| CBM1 | 0 | 0 | 0.011 | 0.0033 | 2191 | 19.9 |
| CBM2 | 0 | 0 | 0 | 0 | 5123 | 33.2 |
| GRD1 | 0.576 | 0.016 | 0.987 | 0.0036 | 469 | 4.1 |
| GRD2 | 0.587 | 0.016 | 0.993 | 0.0026 | 471 | 4.2 |
| GRD3 | 0.062 | 0.0076 | 0.363 | 0.015 | 2300 | 83.5 |
| GRD4 | 0.01 | 0.0031 | 0.083 | 0.0087 | 5365 | 150.5 |

Now we compare GRD and CBM in terms of the quality of estimation. Table 1 gives the estimated mean-squared error (MSE) for each setting based on 1000 independent replications described above. The estimates for GRD were obtained using the methods described earlier in this subsection while the results for CBM were obtained from the simulations performed for Section 4.1.2. It is clear that CBM results in superior estimation. In particular, note that using the setting CBM1 results in better estimates of $E(\mu|y)$ and $E(\lambda|y)$ than using setting GRD4 while using approximately half the average simulation effort [2191 (19.9) versus 5365 (150.5)]; see Table 2.

Consider GRD4 and CBM2. Note that these two settings have comparable average simulation effort. The MSE for $\bar{\mu}_n$ using GRD was 0.000134 (s.e. = $9.2 \times 10^{-6}$) and for CBM we observed an MSE of 0.0000373 ($1.8 \times 10^{-6}$). Now consider $\bar{\lambda}_n$. The MSE based on using GRD was 0.00165 ($1.2 \times 10^{-4}$) and for CBM we observed an MSE of 0.000393 ($1.8 \times 10^{-5}$). Certainly, the more variable simulation effort of GRD contributes to this difference but so does the default use of burn-in.

Recall that we employed a sequential sampler to draw from the target distribution implying that the Markov chain is stationary and hence burn-in is unnecessary. To understand the effect of using burn-in we calculated the estimates of $E(\mu|y)$ using the entire simulation; that is, we did not discard the first $l$ draws of each of the $m$ parallel chains. This yields an estimated MSE of 0.0000709 ($4.8 \times 10^{-6}$) for GRD4. Thus, the estimates using GRD4 still have an estimated MSE 1.9 times larger than that obtained using CBM2. The standard errors of the MSE estimates show that this difference is still significant, indicating CBM, in terms of MSE, is still a superior method for estimating $E(\mu|y)$. Similarly, for estimating $E(\lambda|y)$ the MSE using GRD4 without discarding the first half of each chain is 2.1 higher than that of CBM2.

Toy examples are useful for illustration; however, it is sometimes difficult to know just how much credence the resulting claims should be given. For this reason, we turn our attention to a setting that is "realistic" in the sense that it is similar to the type of setting encountered in practice. Specifically, we do not know the true values of the posterior expectations and implementing a reasonable MCMC strategy is not easy. Moreover, we do not know the convergence rate of the associated Markov chain.

## 5. A HIERARCHICAL MODEL FOR GEOSTATISTICS

We consider a data set on wheat crop flowering dates in the state of North Dakota (Haran et al., 2007). These data consist of experts' model-based estimates for the dates when wheat crops flower at 365 different locations across the state. Let $D$ be the set of $N$ sites and the estimate for the flowering date at site $s$ be $Z(s)$ for $s \in D$. Let $X(s)$ be the latitude for $s \in D$. The flowering dates are generally expected to be later in the year as $X(s)$ increases so we assume that the expected value for $Z(s)$ increases linearly with $X(s)$. The flowering dates are also assumed to be spatially dependent, suggesting the following hierarchical model:

$$Z(s) \mid \beta, \xi(s) = X(s)\beta + \xi(s) \quad \text{for } s \in D,$$
$$\xi \mid \tau^2, \sigma^2, \phi \sim N(0, \Sigma(\tau^2, \sigma^2, \phi)),$$

where $\xi = (\xi(s_1), \ldots, \xi(s_N))^T$ with $\Sigma(\tau^2, \sigma^2, \phi) = \tau^2 I + \sigma^2 H(\phi)$ and $\{H(\phi)\}_{ij} = \exp((-\|s_i - s_j\|)/\phi)$, the exponential correlation function. We complete the specification of the model with priors on $\tau^2$, $\sigma^2$, $\phi$, and $\beta$,

$$\tau^2 \sim \text{IG}(2, 30), \quad \sigma^2 \sim \text{IG}(0.1, 30),$$
$$\phi \sim \text{Log-Unif}(0.6, 6), \quad \pi(\beta) \propto 1.$$

Setting $Z = (Z(s_1), \ldots, Z(s_N))$, inference is based on the posterior distribution $\pi(\tau^2, \sigma^2, \phi, \beta \mid Z)$. Note that MCMC may be required since the integrals required for inference are analytically intractable. Also, samples from this posterior distribution can then be used for prediction at any location $s \in D$.

Consider estimating the posterior expectation of $\tau^2, \sigma^2, \phi$ and $\beta$. Unlike the toy example considered earlier, these expectations are not analytically available. Sampling from the posterior is accomplished via a Metropolis–Hastings sampler with a joint update for the $\tau^2$, $\phi$, $\beta$ via a three-dimensional independent Normal proposal centered at the current state with a variance of 0.3 for each component and a univariate Gibbs update for $\sigma^2$.

To obtain a high-quality approximation to the desired posterior expectations we used a single long run of 500,000 iterations of the sampler and obtained 23.23 (0.0426), 25.82 (0.0200), 2.17 (0.0069) and 4.09 (4.3e-5) as estimates of the posterior expectations of $\tau^2$, $\sigma^2$, $\phi$ and $\beta$, respectively. These are assumed to be the truth. We also recorded the 10th, 30th, 70th and 90th percentiles of this long run for each parameter.

Our goal is to compare the finite-sample properties of GRD and CBM in terms of quality of estimation and overall simulation effort. Consider implementation of GRD. We will produce 100 independent replications using the following procedure. For each replication we used $m = 4$ parallel chains from four different starting values corresponding to the 10th, 30th, 70th and 90th percentiles recorded above. A minimum total simulation effort of 1000 (250 per chain) was required. Also, no burn-in was employed. This is consistent with our finding in the toy example that estimation improved without using burn-in. Each replication continued until $\hat{R}_{0.975} \leq 1.1$ for all of the parameter estimates. Estimates of the posterior expectations were obtained by averaging draws across all four parallel chains.

Now consider the implementation of CBM. For the purpose of easy comparison with GRD, we ran a total of 400 independent replications of our MCMC sampler, where the 10th, 30th, 70th and 90th percentiles of the parameter samples from the long run were used as starting values for 100 replications each. Each replication was simulated for a minimum of 1000 iterations so $p(n) = \varepsilon I(n < 1000)$. Thus the minimum simulation effort is the same as that for GRD. Using (7), a single replication (chain) was terminated when each of the half-widths of a 95% interval was smaller than 0.5, 0.5, 0.05 and 0.05 for the estimates of the posterior expectations of $\tau^2$, $\sigma^2$, $\phi$ and $\beta$, respectively. These thresholds correspond to reasonable desired accuracies for the parameters. If the half-width was not less than the cutoff, then 10 iterations were added to the chain before checking again.

TABLE 3
*Summary of estimated mean-squared error obtained using CBM and GRD for the model of Section 5*

| Method | GRD | | CBM | |
|---|---|---|---|---|
| Parameter | MSE | S.E. | MSE | S.E. |
| $E(\tau^2\|z)$ | 0.201 | 0.0408 | 0.0269 | 0.00185 |
| $E(\sigma^2\|z)$ | 0.0699 | 0.0179 | 0.00561 | 0.00039 |
| $E(\phi\|z)$ | 0.00429 | 0.00061 | 0.000875 | 5.76e-05 |
| $E(\beta\|z)$ | 1.7e-07 | 3.09e-08 | 3.04e-08 | 1.89e-09 |

Standard errors (S.E.) shown for each estimate.

The results from our simulation study are summarized in Table 3. Clearly, the MSE for estimates using GRD are significantly higher than the MSE for estimates obtained using CBM. However, CBM required a greater average simulation effort 31,568.9 (177.73) than did GRD 8,082 (525.7). To study whether the CBM stopping rule delivered confidence intervals at the desired 95% levels, we also estimated the coverage probabilities for the intervals for the posterior expectations of $\tau^2$, $\sigma^2$, $\phi$ and $\beta$, which were 0.948 (0.0112), 0.945 (0.0114), 0.912 (0.0141) and 0.953 (0.0106), respectively. The coverage for all parameters is fairly close to the desired 95%.

Finally, we note that this simulation study was conducted on a Linux cluster using R (Ihaka and Gentleman (1996)), an MCMC package for spatial modeling, spBayes (Finley, Banerjee and Carlin (2007)) and the parallel random number generator package rlecuyer (L'Ecuyer et al., 2002).

## 6. DISCUSSION

In our view, the point of this paper is that those examining the results of MCMC computations are much better off when reliable techniques are used to estimate MCSEs and then *the MCSEs are reported*. An MCSE provides two desirable properties: (1) It gives useful information about the quality of the subsequent estimation and inference; and (2) it provides a theoretically justified, yet easily implemented, approach for determining appropriate stopping rules for their MCMC runs. On the other hand, a claim that a test indicated the sampler "converged" is simply nowhere near enough information to objectively judge the quality of the subsequent estimation and inference. Discarding a set of initial draws does not necessarily improve the situation.

A key requirement for reporting valid Monte Carlo standard errors is that the sampler mixes well. Find-ing a good sampler is likely to be the most challenging part of the recipe we describe. We have given no guidance on this other than one should look within the class of geometrically ergodic Markov chains if at all possible. This is an important matter in any MCMC setting; that is, a Markov chain that converges quickly is key to obtaining effective simulation results in finite time. Thus there is still a great deal of room for creativity and research in improving samplers but there are already many useful methods that can be implemented for difficult problems. For example, one of our favorite techniques is simulated tempering (Geyer and Thompson (1995), Marinari and Parisi (1992)) but many others are possible.

## ACKNOWLEDGMENTS

The authors are grateful to Brian Caffo, Andrew Finley, Andrew Gelman, Charlie Geyer, Jim Hobert, Jim Hodges and Ron Neath for helpful conversations about this paper. The third author's research supported by NSF Grant DMS-08-06178.